\documentclass [a4paper, 11pt] {amsart}
\usepackage [latin1]{inputenc}
\usepackage[all]{xy}
\usepackage{amsmath}
\usepackage{amsfonts}
\usepackage{latexsym}  
\usepackage{amssymb}
\usepackage{graphics}
\usepackage{graphicx}
\newtheorem{teorema}{Theorem}[section]
\newtheorem{Lemma}[teorema]{Lemma}
\newtheorem{propos}[teorema]{Proposition}
\newtheorem{corol}[teorema]{Corollary}
\newtheorem{ex}{Example}[section]
\newtheorem{rem}{Remark}[section]
\newtheorem{defin}[teorema]{Definition}
\def\bt{\begin{teorema}}
\def\et{\end{teorema}}
\def\bp{\begin{propos}}
\def\ep{\end{propos}}
\def\bl{\begin{Lemma}}
\def\el{\end{Lemma}}
\def\bc{\begin{corol}}
\def\ec{\end{corol}}
\def\br{\begin{rem}\rm}
\def\er{\end{rem}}
\def\bex{\begin{ex}\rm}
\def\eex{\end{ex}}
\def\bd{\begin{defin}}
\def\ed{\end{defin}}

\def\R{{\mathbb R}}    
\def\N{{\mathbb N}}      
\def\C{{\mathbb C}}

\def\D{\Delta}

\begin {document}
\title{Discrete sequences in unbounded domains}
\author[Alberto Saracco]{Alberto Saracco}
               \address{Dipartimento di Matematica e Informatica, Universit\`a di Parma, Parco Area delle Scienze 53/A, I-43124 Parma, Italy}\email{alberto.saracco@unipr.it}
                             
               \date{\today}

\begin{abstract} Discrete sequences with respect to the Kobayashi distance in a strongly pseudoconvex bounded domain $D$ are related to Carleson measures by a formula that uses the Euclidean distance from the boundary of $D$.

Thus the speed of escape at the boundary of such sequence has been studied in details for strongly pseudoconvex bounded domain $D$.

In this note we show that such estimations completely fail if the domain is not bounded.\end{abstract}
 \maketitle

Let $(X,d)$ be a metric space.  A sequence $\Gamma=\{x_j\}\subset
X$ of points in $X$ is {\sl uniformly discrete} if there exists $\delta>0$
such that $d(x_j,x_k)\ge\delta$ for all $j\ne k$. 

In this case $\inf\limits_{j\ne k} d(x_j,x_k)$ is the {\sl separation constant} of~$\Gamma$.\vspace{0.3cm}

Finite unions of uniformly discrete sequences in a strongly pseudoconvex bounded domain $D\subset\C^n$ with respect to the Kobayashi distance are a mean to produce explicit examples of discrete Carleson measures.

\begin{teorema}\label{0}Let $D\Subset{\C}^n$ be a strongly pseudoconvex bounded  domain, and let $\Gamma=\{z_j\}_{j\in\N}$ be a sequence
in~$D$. Then $\Gamma$ is a finite union of uniformly discrete (with respect to the Kobayashi
distance) sequences if and only if 
$$ \sum\limits_{z_j\in\Gamma} d(z_j,\partial  D)^{n+1}\delta_{z_j}$$
is a Carleson measure of $A^p(D)$, where $\delta_{z_j}$ is the Dirac measure at~$z_j$
and $d(\cdot,\partial  D)$ is the Euclidean distance from the boundary.
\end{teorema}

This theorem, by Abate and the author~\cite{AS} is a generalization of results by Massaneda~\cite{M}, Jevti\'c, Massaneda and Thomas~\cite{JMT} and Duren and Weir~\cite{DW}. These results motivated the study of the boundary behaviour, with respect to the Euclidean distance, of uniformly discrete sequences (with respect to the Kobayashi distance) in a strongly pseudoconvex bounded domain $D\subset\mathbb C^n$.  In~\cite{DSV},~\cite{DW},~\cite{Mac} and~\cite{AS} were obtained results in this direction for $D$ being the unit ball or a strongly pseudoconvex bounded domain or even a hyperbolic domain with finite Euclidean volume.\vspace{0.3cm}

We cite here two theorems proved by Abate and the author in~\cite{AS}. The first one generalizes the results obtained by Duren, MacCluer, Schuster, Vukoti\'c and Weir for the unit ball, to strongly pseudoconvex bounded domains (and is a sharp estimation). The second one is valid in a more general setting but provides a worst estimation of the boundary behaviour of uniformly discrete sequences.

\begin{teorema}\label{1} Let $D\Subset{\C} ^n$ be a strongly pseudoconvex bounded domain.

Let $\Gamma=\{z_j\}\subset D$ be an uniformly discrete sequence (with respect to the Kobayashi distance) with $d(z_j,\partial  D)<1$ for all~$z_j\in\Gamma$. Then
\[
\sum_{z_j\in\Gamma} d(z_j,\partial  D)^{n}\, h\left(-{1\over\log
d(z_j,\partial  D)}\right)<+\infty
\]
for any increasing function $h\colon\R^+\to{\R} ^+$ such that
\[
\sum_{m=1}^{+\infty}h\left({1\over m}\right)<+\infty\;.
\]
\end{teorema}

\begin{teorema}\label{2} Let $D\subset{\C} ^n$ be a hyperbolic domain with finite Euclidean
volume.

Let $\Gamma=\{z_j\}\subset D$ be an uniformly discrete sequence  (with respect to the Kobayashi distance) with $d(z_j,\partial  D)<1$ for all~$z_j\in\Gamma$. Then
\[
\sum_{z_j\in\Gamma} d(z_j,\partial  D)^{2n}\, h\left(-{1\over
\log d(z_j,\partial  D)}\right)<+\infty
\]
for any increasing function $h\colon{\R} ^+\to{\R} ^+$ such that
\[
\sum_{m=1}^{+\infty}h\left({1\over m}\right)<+\infty\;.
\]
\end{teorema}

One might hope that the boundedness (or limited volume) condition is unnecessary, at least in cases where the hyperbolicity of the domain imply some form of boundedness. For example, (complete) hyperbolic convex (or even $\C$-convex) domains are known to be biholomorphically equivalent to a bounded domain (see~\cite{BS} for the convex case,~\cite{NPZ} and~\cite{NS} for the $\C$-convex case).

If such results were true this would set a much general result for any convex (respectively $\C$-convex) domain, since such a domain $D\subset\C^n$ is linearly equivalent to the product of $\C^{n-k}$ and a (complete) hyperbolic convex  (respectively $\C$-convex) domain $D'\subset\C^k$ and the Kobayashi pseudodistance in $D$ would coincide with the Kobayashi distance on the hyperbolic factor $D'$ (see again~\cite{BS},~\cite{NPZ} and~\cite{NS}).

Unfortunately, this is actually not possible. We give an example showing that for the half plane $H\subset\C$ and then generalize the result to a huge class of unbounded domains.

\section{An example}
Let $\mathbb D\subset \C$ be the unit disc. $\mathbb D$ is biholomorphic to the right half plane $H=\{z\in \C\ |\ \mathfrak{Re}\,z>0\}$ via the function
$$\varphi(z)\ =\ \frac{1-z}{1+z}\,.$$

Let us consider the set of points at a fixed Euclidean distance $\varepsilon>0$ from the boundary of $H$:
$$l_\delta=\{\varepsilon+iy\ |\ y\in\R\}\,.$$
Its image in $\D$ via $\varphi^{-1}$ is an ellipse $E_\varepsilon$ internally tangent to $\partial \D$ of center
$$C\ =\ -\frac{\varepsilon}{1+\varepsilon}\,.$$
Let us fix a $\delta>0$. We define a discrete sequence of points $\Gamma=\{z_j\}_{j\in\N}$ in $\D$ as follows. Fix $z_0\in \D\cap E_\varepsilon$. Since $E_\varepsilon$ is tangent to $\partial \D$, for every $k\in \mathbb N$ the set $E_{k,\varepsilon}=E_\varepsilon\cap\{z\in \D \ |\ d(z_0,z)=k\delta\}$ is not empty. Hence we can chose $z_k\in E_{k,\varepsilon}$ and by the triangular inequality $d(z_h,z_k)\geq\delta$ for every $h\neq k\in \mathbb N$, i.e.\ $\Gamma=\{z_j\}_{j\in\mathbb N}$ is a uniformly discrete sequence with separation constant $\delta$. The points $z_j$ approach the boundary of $\D$, as $j$ goes to infinite.

If we look at the very same sequence in $H$, i.e.\ if we define $w_j=\varphi^{-1}(z_j)$, then $\varphi^{-1}(\Gamma)=\{w_j\}_{j\in\mathbb N}$ is a uniformly discrete sequence in $H$ with separation constant $\delta$, since $\varphi:H\to \D$ is a biholomorphism, and $d(w_j,\partial H)=\varepsilon$ for each $j\in\mathbb N$.

Hence if $$F\colon{\R} ^+\to{\R} ^+$$ is any function such that $F(x)>0$ for every $x>0$, then
$$\sum_{w_j\in\varphi^{-1}(\Gamma)} F(d(w_j,\partial H))\ =\ \sum_{j=0}^\infty F(\varepsilon)\ =\ +\infty\,.$$

In particular a result such as Theorem~\ref{1} or~\ref{2}, valid in $\mathbb D$, fails in $H$, despite the two domains being biholomorphic. This is obviously because Euclidean distance from boundary is not preserved by biholorrphisms.

\section{A theorem}

The argument of the previous example can be generalized to prove the following

\begin{teorema} Let $D\subset\mathbb C^n$ be a Kobayashi hyperbolic domain such that there exists $c>0$ for which
$D_c\ =\ \{z\in D\ | \ d(z,\partial D)\geq c\}\subset D$ is non compact. Then for any $\delta>0$ there exists a uniformly discrete sequence $\Gamma=\{z_j\}_{j\in\N}$ in $D$ with separation constant $\delta$ such that
$\Gamma\subset D_c$.
In particular
$$\sum_{w_j\in\varphi^{-1}(\Gamma)} F(d(w_j,\partial D))\ =\ +\infty\,,$$
for every increasing function $F\colon{\R} ^+\to{\R} ^+$ such that $F(x)>0$ for every $x>0$.
\end{teorema}

\br Thus no equivalent of theorems~\ref{1} or~\ref{2} holds for such a domain.
\er

\br Note that if the set $D_c$ is non compact, i.e.\ unbounded since it is closed, then $D\setminus D_c$ has not a finite Euclidean volume, hence the hypothesis of theorem~\ref{2} does not hold.
\er

\begin{proof}Let $c>0$ be as in the hypothesis. Fix $\delta>0$. We will construct the required uniformly discrete sequence $\Gamma$ (with separation constant $\delta$) by induction.

Since $D_c$ is non compact, it is not empty. Fix $z_0\in D_c$.

Let $k\in\N$. Suppose for all $j<k$ the points $z_j\in D_c$ has been choosen in such a way that any two such points are at least $\delta$-apart from each other with respect to the Kobayashi distance. Since $D$ is Kobayashi hyperbolic, for every $j<k$ $B_j=B_\delta(z_j)$, the Kobayashi ball of radius $\delta$ and center $z_j$ is relatively compact in $D$. Hence $\cup_j B_j$ is relatively compact in $D$ and thus $D_c\setminus\cup_j B_j$ is not empty. By choosing $z_k\in D_c\setminus\cup_j B_j$, $z_k$ is at least $\delta$-apart from each other $z_j$, $j<k$, and the inductive step is done.

Hence the claimed $\Gamma$ is constructed.\end{proof}


\begin{thebibliography}{}

\bibitem{AS} Marco Abate \and Alberto Saracco, \textit{Carleson measures and uniformly discrete sequences in strongly pseudoconvex domains}, J.\ London Math.\ Soc.\ \textbf{84} (2011), 587--605.

\bibitem{BS} Filippo Bracci \and Alberto Saracco, \textit{Hyperbolicity in unbounded convex domains}, Forum Math.\ \textbf{5} (2009), 815--826.

\bibitem{DSV} Peter Duren, Alexander Schuster \and Dragan Vukoti\'c, \textit{On uniformly discrete sequences in the disk}, in
\textsl{Quadrature domains and applications,} eds. P. Ebenfeld, B. Gustafsson, D. Khavinson and M. Putinar, Birkh\"auser, Basel, 2005, pp. 131--150.

\bibitem{DW} Peter Duren \and Rachel Weir, \textit{The pseudohyperbolic metric and Bergman spaces in the ball}, Trans.\ Amer.\ Math.\ Soc.\ \textbf{359} 1 (2007), 63--76.

\bibitem{JMT} Miroljub Jevti\'c, Xavier Massaneda \and Pascal J. Thomas, `Interpolating sequences for weighted Bergman 
spaces of the ball', \textsl{Michigan Math. J.} 43 (1996), 495--517.

\bibitem{Mac} Barbara D.\ MacCluer, \textit{Uniformly discrete sequences in the ball}, J.\ Math.\ Anal.\ Appl.\ \textbf{318} (2006), 37--42.

\bibitem{M} Xavier Massaneda, `$A^{-p}$ interpolation in the unit ball', \textsl{J. London Math. Soc.} 52 (1995),
391--401.

\bibitem{NPZ} Nikolai Nikolov, Peter Pflug \and W\l odzimierz Zwonek, \textit{An example of a bounded ${\mathbb C}$-con\-vex domain which is not biholomorphic to a convex domain}, Math.\ Scand.\ \textbf{102} 1 (2008), 149--155.

\bibitem{NS} Nikolai Nikolov \and Alberto Saracco, \textit{Hyperbolicity of unbounded $\C$-convex domains}, C.\ R.\ Acad.\ Bulgare Sci.\ \textbf{60} 9 (2007), 737--748.


\end{thebibliography}
\end{document}